\documentstyle{amsart}
\def\blacksquare{\quad \vrule height 6pt width 6pt}
\newenvironment{prf}{{\bf\underline{proof}}:}{{\samepage\hfill
\hbox{\blacksquare}}\par}

\newtheorem{prp}{Proposition}
\newtheorem{thm}{Theorem}

\def\k{\mbox{\bf k}}				%letters

\def\g{\mbox{$\frak g$}}
\def\m{\mbox{$\frak m$}}

\def\Z{\mbox{$\Bbb Z$}}

\def\zt{\hbox{$\Z_2$}}				%Z-two
\def\ztz{\hbox{$\Z_2\times\Z$}}			%Z-two tensor Z

\def\bdeg{\operatorname{bid}}				%operators

\def\coder{\operatorname{Coder}}
\def\im{\operatorname{im}}
\def\sh{\operatorname{Sh}}				%"shuffle"
\def\hom{\mbox{\rm Hom}}      			%no space after
\def\ker{\mbox{\rm Ker}\, }   			%space after
\def\span{\mbox{\rm span}}    			%no space after

\def\linf{\mbox{$L_\infty$}}
\def\ainf{\hbox{$A_\infty$}}
\def\and{\hbox{ and }}

\def\ph{\varphi}					%the variable phi
\def\ra{\rightarrow}				% -->
\def\dirlim{lim \kern-18pt\lower1ex\hbox{$\longleftarrow$}\ } %direct limit
\def\iso{\kern.35em{\raise3pt\hbox{$\sim$}\kern-1.1em\to}\kern.3em}  %isomorphic
\def\of{\circ}                			%function composition
\def\tns{\otimes}             			%regular tensor product
\def\Tns{\bigotimes}				%big tensor product

\def\mcom{,\cdots,}					%         ,...,
\def\mplus{+\cdots+}				%         +...+
\def\mwedge{\wedge\cdots\wedge}			%     wedge...wedge

\def\leftv{v_1\mcom v_j}				%subscripts of v
\def\midv{v_{j+1}\mcom v_{j+l}}
\def\rightv{v_{j+l+1}\mcom v_{n}}

\def\sigsign{\s{\sigma}\epsilon(\sigma)}		%signs notation
\def\e#1{\mbox{$|#1|$}}
\def\s#1{\mbox{$(-1)^{#1}$}}

\def\br#1#2{\left\{#1,#2\right\}}			%brackets
\def\ip#1#2{\left<#1,#2\right>}			%inner product notation

\def\permarg#1#2#3{\mbox{$ {#1}_{\sigma(#2)}, \cdots, {#1}_{\sigma(#3)}$}}  %permuted arguments
\begin{document}
\nocite{mar,sta3,getz2,lod,conn,kast,seib,hoch,umb,aksz}
\author{Michael Penkava}
\address{University of Winsconsin-Eau Claire\\
Eau Claire, WI 54702-4004}
\email{penkavmn@@uwec.edu}
\author{Lynelle Weldon}
\address{Andrews University\\
Berrien Springs, MI 49104-0350}
\email{weldon@@andrews.edu}
\title{Infinity Algebras, Massey Products, and Deformations}
\maketitle
\section{Introduction}
The notion of a Massey F-product was introduced in \cite{fl} in order
to give a unified description of various deformation problems in terms
of an extension of the usual Massey product. The usual Massey product which
arises in algebra deformation theory describes the conditions under which
an infinitesimal deformation can be prolonged to a higher order or formal
deformation. The Massey F-product was introduced in order to solve some
more general deformation problems. Associated to each Massey F-product is 
a \zt-graded coalgebra, which depends on the problem being described.
In this paper we consider some examples
which arise in connection with the problem of deforming \ainf\ and \linf\
algebras, as well as the more restricted problem of deforming associative
algebras into \ainf\ algebras and Lie algebras into \linf\ algebras. 
The general example of deforming \ainf\ and \linf\ algebras turns out to
give the same coalgebra structure as arises in the problem of prolonging an
infinitesimal deformation of
a \zt-graded Lie algebra into a higher order deformation. The restricted
problem yields a more interesting coalgebra structure, which we will
describe here. 

We also consider deformations of \ainf\ and \linf\ algebras with a base 
given by a \zt-graded commutative algebra, and prove an extension of 
a theorem in \cite{fl}, that deformations of Lie algebras
with a base are determined by triviality of Massey F-products.  

In order to make the reading of this paper more self-contained, we include
definitions of \ainf\ and \linf\ algebras. \ainf\ algebras were first
introduced by J. Stasheff
in \cite{sta1,sta2}, while the notion of an \linf\ algebra
first appeared in \cite{ss}. Our definitions are based on the exposition
in \cite{pen1}, but the articles \cite{ls,lm,mar,getz} contain
an introduction to infinity algebras as well.

\section{Notation}
In the following we shall be considering vector spaces graded by some
group equipped  with a \zt-valued inner product.
For any graded vector space $V$, define maps $S:V\Tns V\ra V\Tns V$
and $C:V\Tns V\Tns V\ra V\Tns V\Tns V$ by
$S(u\tns v)=\s{uv}(v\tns u)$ and $C(u\tns v\tns w)=\s{u(v+w)}
v\tns w\tns u$, where we adopt the convention that expressions like
$\s{uv}$ stand for the sign determined by the inner product of the
gradings of $u$ and $v$. These signs appear, for example, 
when considering tensor products of mappings from  graded spaces.

A graded coassociative coalgebra is a graded vector space $V$
equipped with 
a degree zero mapping (comultiplication) $\Delta:V\ra V\Tns V$
satisfying $(1\tns\Delta)\circ\Delta=(\Delta\tns1)\circ\Delta$.
It is cocommutative if $S\circ\Delta=\Delta$.
In what follows, we shall consider \zt-graded vector spaces, but
everything said applies to the \Z-graded case as well. If $V$ is
a \zt-graded vector space, then the tensor and exterior (co)algebras
are \ztz-bigraded algebras. 
We shall encounter two different inner products on \ztz,
the {\em usual inner product} given by 
$\ip{(\bar k,m)}{(\bar l,n)}=\overline{kl+mn}$, and the inner product
$\ip{(\bar k,m)}{(\bar l,n)}=\overline{(k+m)(l+n)}$, which we will 
call the {\em good inner product}, for reasons which will appear later.
(In the formulas above, the bar stands for integers $\mod 2$.)
We consider the two inner products as giving different gradings on $T(V)$.
In both gradings, a homogeneous element $v$ of bidegree $(\bar k,m)$
will be called odd if $k+m$ is an odd integer.
If $v$ is an element of a \ztz-graded space $V$ then
$(-1)^v=(-1)^{|v|}$ is the sign determined by the parity of $v$.

\section{\ainf\ Algebras}
Suppose $V$ is a \zt-graded vector space, and consider the tensor
coalgebra $T(V)=\bigoplus_{n=1}^\infty V^n$. $T(V)$ has a natural
\ztz-bigrading, the bidegree of $v\in V^n$ being given by
$\bdeg(v)=(\e v,n)$. 
An \ainf\ algebra is a vector space $V$ with an odd map $\mu\in\hom(T(V),V)$,
in other words a sequence of maps $\mu_k:V^k\ra V$, which satisfy
for $n\ge 1$,
\begin{equation}\label{ainfeq}
\sum_{k+l=n+1 \atop 0\leq j \leq k-1}
\s r
\mu_k(\leftv,\mu_l(\midv),\rightv)=0,
\end{equation}
where $r=l(v_1\mplus v_j)+j(l-1)+(k-1)l$.
An \ainf\ algebra may be seen as a generalization of an associative
algebra, and in particular, if $\mu$ consists only of an $\mu_2$ part,
then it is an associative algebra. To explain the origin of the
signs above consider the complex, $C(V)=\hom(T(V),V)$ with the
usual \ztz-grading. 
The bigrading on $C(V)$ is given by $\bdeg(\ph)=(\e\ph,k-1)$ for 
$\ph\in C^k(V)=\hom(V^k,V)$.
There is a natural isomorphism $C(V)\cong \coder(T(V))$,
the \ztz-graded coderivations of $T(V)$.
Since $\coder(T(V))$ is a
graded Lie algebra with respect to the usual \ztz-grading, 
$C(V)$ has a Lie bracket, which is given by
\begin{multline}\label{tensbra}
[\ph, \psi](v_1,\cdots,v_n)=\\
\sum_{0\leq j \leq k-1}
\s{r}
 \ph(\leftv,\psi(\midv),\rightv)\\
-(-1)^{\ph\psi + (k-1)(l-1)}\sum_{0\leq j\leq l-1}\!\!\!
\s{s}
\psi(\leftv,\ph(\midv),\rightv).
\end{multline}
for $\ph\in C^k(V)$, 
$\psi\in C^l(V)$, where $n=k+l-1$, $r=\psi(v_1+\cdots v_{j})+j(l-1)$
and $s=\ph(v_1+\cdots+v_{j}) +j(k-1)$.
This bracket is called the Gerstenhaber bracket. 
Lie algebras with
the usual \ztz-grading have some undesirable properties, but there is
a procedure for modifying the bracket which transforms a Lie algebra
with the usual \ztz-grading into a Lie algebra with the good \ztz-grading.
The modified Gerstenhaber bracket is defined by
$$\{\ph,\psi\}=(-1)^{(k-1)\psi}[\ph,\psi].$$
An \ainf\ structure is determined by an odd element $\mu$ satisfying
$\br \mu\mu=0$. This definition is equivalent to the relations given by
equation (\ref{ainfeq}).
Then $\delta(\ph)=\{\mu,\ph\}$ defines a
differential on $\hom(T(V),V)$ and gives it a DGLA 
(Differential Graded Lie Algebra) structure.
The property that $\delta^2=0$ follows from the Jacobi identity. This
property is what makes the good grading behave better than the usual
one. 
We shall call $\mu$ the codifferential determining the \ainf\
structure on $V$, although it is not
precisely a codifferential, since it is given by the vanishing of 
the modified bracket, rather than the bracket of coderivations.

\section{\linf\ Algebras}
To define an \linf\ algebra, we begin with the exterior coalgebra
$\bigwedge V$ of a \zt-graded vector space $V$. 
With respect to the usual \ztz-grading, $\bigwedge V$ is a cocommutative
coalgebra.
For $v_1\mwedge v_n \in \bigwedge V$ and any permutation $\sigma$,
$\epsilon (\sigma ; v_1, \cdots, v_n)$ is defined by the equation
\begin{equation}
 v_1\mwedge v_n = (-1)^\sigma \epsilon(\sigma ; v_1, \cdots, v_n)
  v_{\sigma(1)}\mwedge v_{\sigma(n)},
\label{eps}
\end{equation}
where $(-1)^\sigma$ is the sign of the permutation $\sigma$.  
For ease of notation, $\epsilon(\sigma;v_1,\cdots,v_n)$ 
will be denoted by $\epsilon(\sigma)$.
An \linf\ algebra is given by an odd
map $l\in C(V)=\hom(\bigwedge V,V)$, that is, 
with a sequence of maps
$l_k:\bigwedge^k V \to V$, with $\e {l_k}=k$,
which satisfy
\begin{equation}\label{linfdef}
\sum _{{k+l =n+1}\atop{\sigma\in\sh(k,n-k)}}
\s{\sigma}\epsilon(\sigma)\s{(k-1)l}
l_l(l_k(\permarg{v}{1}{k}),\permarg{v}{k+1}{n})=0
\end{equation} 
for $n\ge1$.  
The set of relations above will be referred to as the generalized
Jacobi identity.
This structure is also called a homotopy Lie
algebra and may be viewed as a generalization of a Lie algebra as
follows.  If $l_k=0$ for $k>2$ then $V$ has a DGLA
structure with differential $l_1$ and bracket given by $l_2$.  
This generalization was motivated by deformation theory \cite{ss} and
has applications in quantum mechanics and string theory \cite{ber,WitZwie,Z}.
 
The bigrading on $C(V)$ is given by 
$\bdeg(\ph)=(\e\ph,k-1)$ for 
$\ph\in C^k(V)=\hom(V^k,V)$.
There is a natural isomorphism $C(V)\cong \coder(\bigwedge V)$,
equipping it with a \ztz-graded Lie bracket defined  by
\begin{multline}
[\ph,\psi](v_1,\cdots,v_n)=\\
\sum_{\sigma\in\sh( l,k-1)}\sigsign
 \ph(\psi(\permarg{v}{1}{l}),\permarg{v}{l+1}{n})\\
-(-1)^{\psi\ph+(k-1)(l-1)}\!\!\!\sum_{\sigma\in\sh(k,l-1)}\sigsign
 \psi\ph(\permarg{v}{1}{k}),\permarg{v}{k+1}{n}),
\end{multline}
$\ph\in\hom(\bigwedge^k(V),V),\psi\in\hom(\bigwedge^lV,V),$
and\ $v_1,\cdots,
v_n\in V$, where $k+l=n+1$.
This bracket is simply the bracket of $\ph$ and $\psi$ as 
coderivations with respect to the usual \ztz-grading.
As in the \ainf\ algebra case, a modified  bracket is defined by
\begin{equation}
\{\ph,\psi\}=(-1)^{(k-1)\psi}[\ph,\psi].
\label{modbrack}
\end{equation}
The modified bracket gives $C(V)$ the structure of a \ztz-graded
Lie Algebra with respect to the good inner product on \ztz.
Then equation (\ref{linfdef}) is simply the property that  
$\br ll=0$. As before, we define a differential on $C(V)$ by
$\delta(\ph)=\br l\phi$. 
Thus $C(V)$ inherits
the structure of a DGLA.
We shall call $l$ the codifferential determining the \linf\ structure
on $V$. In this fashion, we say that \ainf\ and \linf\ structures
are determined by codifferentials, although this is not precisely true.
\section{\label{SSdirect}Definition of Massey F-Products}
Let $L$ be a DGLA with commutator 
 $\mu\colon L\tns L\to L,\ \mu(u,v)=[u,v]$, and
differential  $\delta\colon L\to L$.  Denote the
cohomology of $L$ with respect to $\delta$ by $H=\bigoplus_iH^i$.
Let $F$ be a graded cocommutative coassociative coalgebra with
comultiplication operator  $\Delta\colon F\to F\tns F$. 
Suppose also that a filtration  $F_0\subset F_1\subset F$ of $F$ is given, 
such that $F_0\subset\ker \Delta$ and $\im \Delta\subset F_1\tns
F_1$. 
Let $a\colon F_0\to H$, $b\colon F/F_1\to H$ be linear 
maps of degree 1. We say that
$b$ is contained in the Massey $F$-product of $a$
if there exists a degree 1 linear mapping 
$\alpha\colon F_1\to L$ satisfying the condition 
\begin{equation}
 \delta\of\alpha=\mu\of(\alpha\tns\alpha) \of\Delta, 
\label{theeqn}
\end{equation}
and such that the diagrams 
 
\begin{picture}(300,70)
\put(30,5){\begin{picture}(100,70)
 \multiput(2,0)(0,46){2}{\makebox(0,0)[u]{$F_0$}}
 \multiput(0,13)(3,0){2}{\line(0,1){22}}
 \put(15,42){$@>{\ \ \alpha |_{F_0}\ }>>$}
 \put(15,-2){$@>{\ \ \ a\ \ \ }>>$} 
 \put(70,46){\makebox(0,0)[u]{Ker $\delta$}}
 \put(67.3,20){$\Big\downarrow$}
 \put(75,23){$\scriptstyle\pi$}
 \put(90,25){,}
 \put(70,0){\makebox(0,0){$H$}}	
\end{picture}}

\put(200,5){\begin{picture}(100,70)
 \put(0,46){\makebox(0,0)[u]{$F$}}
 \put(0,0){\makebox(0,0)[u]{$F/F_1$}}
 \put(100,0){\makebox(0,0)[u]{$H$}}
 \put(100,46){\makebox(0,0)[u]{Ker $\delta$}}
 \put(17,42){$@>{\ \mu\of(\alpha\tns\alpha)\of\Delta\ }>>$}
 \put(17,-2){$@>{\ \ \ \ \ \ b\ \ \ \ \ \ \ }>>$}
 \put(-3.3,20){$\Big\downarrow$}
 \put(4.4,23){$\scriptstyle\pi$} 
 \put(97.1,20){$\Big\downarrow$}
 \put(104.8,23){$\scriptstyle\pi$}
\end{picture}}

\end{picture}

\smallskip
\noindent
are commutative, where the vertical maps labeled $\pi$ denote the 
projections of each space onto the quotient space. 

Note that from (\ref{theeqn}) it follows that 
$\alpha(F_0)\subset\ker\delta$  and 
$\mu\of(\alpha\tns\alpha)\of\Delta(F)\subset\ker \delta$, 
making upper horizontal maps of the diagrams well-defined (see \cite{fl}).

In the case that $F_1=F$, $b$ does not need to be specified and $a$
is said to satisfy
the condition of triviality of Massey F-products if such an 
$\alpha$ exists.

\section{Graded Lie Algebra Deformations}
Massey F-Products arise in the deformation theory of graded Lie
algebras in the following fashion.  
Let $V$ be a \zt-graded Lie algebra over \k\ 
with bracket $l(v_1,v_2)=[v_1 ,v_2 ].$  
A {\it formal deformation} of $V$ may be defined as a power
series
\begin{equation}
 [g,h]_t = [g,h]
           +\sum_{p=1}^\infty t^p(\gamma_p(g,h)+\theta\beta_p(g,h)),
\label{grdef}
\end{equation}
where $t$ is an even parameter and $\theta$ is an odd parameter.  
In order for the new bracket to be bilinear and antisymmetric 
each $\gamma_i$ and $\beta_i$ must be bilinear and antisymmetric, 
thus  $\gamma_i$ and $\beta_i$ are cochains in the cohomology of $V$.
The bracket of cochains, $\ph\in C^k(V),\psi\in C^l(V),$ is given by
 \begin{multline}
[ \ph,\psi](v_1,\cdots,v_n)=\\
\sum_{\sigma\in\sh(l,k-1)}\sigsign
   \varphi(\psi(\permarg{v}{1}{l}),\permarg{v}{l+1}{n})\\
-(-1)^r
\sum_{\sigma\in\sh(k,l-1)}\sigsign
\psi(\varphi(\permarg{v}{1}{k}),\permarg{v}{k+1}{n}),
\end{multline}
where $r=\varphi\psi+(k-1)(l-1)$ and $k+l=n+1$.
Defining a differential on $C(V)$ by $ \delta\ph=[l,\ph]$
for any $\ph\in C(V)$ 
gives a DGLA structure on $C(V).$
Using this structure, the deformed bracket satisfies the graded Jacobi
identity precisely when
\begin{align}
 \delta\gamma_p&=-{1\over2}\sum_{i=1}^{p-1}
[\gamma_i,\gamma_{p-i}]
\label{gamma} \\
\delta\beta_p&=-{1\over2}\sum_{i=1}^{p-1}[\beta_i,\gamma_{p-i}]
                                        +[\gamma_i,\beta_{p-i}].
\label{beta}
\end{align}
For $p=1,$ the sums in (\ref{gamma}) and (\ref{beta}) 
are empty so $\gamma_1$ and $\beta_1$ are
cocycles representing certain cohomology classes.  The question of
interest is, when does a deformation (\ref{grdef}) exist with
$\gamma_1\in c, \beta_1\in b$ representing two given cohomology
classes, $c$ and $b$?  Such a deformation exists when $a$ satisfies
the condition of triviality of Massey F-products where
\begin{gather*}
 F=\span\{e^1,e^2,\cdots;f^1,f^2,\cdots\},\\
 F_0=\span\{e^1,f^1\},\\
 F_1=F,\\
 \bdeg e^i = (0,0), \bdeg f^i = (1,0),\\
 \Delta e^k = -{1\over2}\sum_{i=1}^{k-i}e^i\tns e^{k-i}, \\
 \Delta f^k = -{1\over2}\sum_{i=1}^{k-i}f^i\tns e^{k-i}
                                     +e^i\tns f^{k-i},\\
a(e^1)=c, a(f^1)=b.
\end{gather*}
Let $G$ be the dual algebra to $F$, with basis $e_k$, $f_k$.
Then $G$ is a commutative algebra with structure given by
multiplication
$$e_ke_l=-{1\over2}e_{k+l},\quad e_k f_l=-{1\over2} f_{k+l},
\quad  f_k f_l=0.$$
If we let $t=e_1$ and $\theta=f_1$   
then $e_k=(-2)^{k-1}t^k,f_l=(-2)^{k-2}t^{k-2}\theta,$ 
and $\theta^2=0.$ 
In this form, $G$ is the maximal ideal in $\k[[t,\theta]]/{\theta^2}$, 
and $F_0^\ast =G/G^2.$
\section{Deformations of \ainf\ and \linf\ Algebras}
A formal deformation of an \linf\ or \ainf\ algebra $V$ with 
codifferential $d$ is given by 
\begin{equation}
\label{infdef}
d_t=\sum_{i=0}^\infty t^i(\gamma_i+\theta\beta_i)
\end{equation}
with the property that $\{d_t,d_t\}=0$.  
Here $\gamma_i, \beta_i$ are cochains in $C(V)$,
$\theta$ is an odd parameter, 
$\alpha_0=d$ denotes the original codifferential determining the 
structure on 
$V$, and 
$\beta_0=0$.  
The property $\{d_t,d_t\}=0$ holds exactly when, 
for all $p$, 
\begin{align}
\delta\gamma_p&=-{1\over2}\sum_{i=1}^{p-1}
\left\{ \gamma_i,\gamma_{p-i}\right\}
\label{lin1}\\
\delta\beta_p&={1\over2}\sum_{i=1}^{p-1} 
\left\{\beta_i,\gamma_{p-i}\right\}-\left\{\gamma_i,\beta_{p-i}\right\}.
\label{lin2}
\end{align}
This gives us the same coalgebra structure as the graded Lie algebra
case (the apparent sign differences work out because of the change in
the bracket being used). 
\section{Lie and Associative Algebra Case} 
The cochains $\gamma_i$ and $\beta_i$ may be written as infinite sums
$\gamma_i=\sum_{m=1}^\infty \ph_{i,m}$ 
and $\beta_i=\sum_{m=1}^\infty \psi_{i,m}$,
where $\bdeg (\ph_{i,m}) = (m,m-1)$, and $\bdeg(\psi_{i,m})=(m-1,m-1)$.
In the case of deforming a Lie algebra into an \linf\ algebra or an
associative algebra into an \ainf\ algebra
$\ph_{0,j}=0$ for $j\neq 2$,
so equations  (\ref{lin1}) and (\ref{lin2}) 
may be broken down into conditions on the $\ph$ and $\psi$ as follows:
\begin{align*}
 \delta\ph_{p,q}&=-{1\over2}\sum_{i=1}^{p-1}\sum_{j=1}^{q+1}
   \left\{\ph_{i,j},\ph_{p-i,q+2-j} \right\} \\
  \delta\psi_{p,q}&={1\over2}\sum_{i=1}^{p-1}\sum_{j=1}^{q+1}
   \left\{\psi_{i,j},\ph_{p-i,q+2-j}\right\}
   -\left\{\ph_{i,j},\psi_{p-i,q+2-j}\right\}.
\end{align*}
Given sequences $c_i$ and $b_i$ of cohomology classes a deformation
(\ref{infdef})  with 
$\ph_{1,i}\in c_i$ and
$\psi_{1,i}\in b_i$ exists if and only if $a$ satisfies
the condition of triviality for the following Massey F-product:
\begin{gather*}
 F = \span\left\{e^{i,j},f^{i,j}\right\}_{i\ge 1, i+j \ge 2}\\
 F_0=\span\left\{e^{i,j},f^{i,j}\right\}_{i\le 1}\\
 F_1=F\\
 \bdeg e^{i,j}=(j,j-2),\bdeg f^{i,j}=(j-1,j-2)\\
 \Delta e^{p,q}=-{1\over2}\sum_{i=1}^{p-1}\sum_{j=1}^{q+1}
 e^{i,j}\tns e^{p-i,q+2-j}\\ 
 \Delta f^{p,q}=-{1\over2}\sum_{i=1}^{p-1}\sum_{j=1}^{q+1}
  f^{i,j}\tns e^{p-i,q+2-j}+e^{i,j}\tns f^{p-i,q+2-j}\\
a(e^{i,j})=
\begin{cases}
   c_j&   i=1\\
   0&     i<1
\end{cases}\\
a(f^{i,j})=
\begin{cases}
  b_j&   i=1\\
  0&     i<1
\end{cases}.
\end{gather*}

The dual algebra $G$ is spanned by $\{e_{i,j},f_{i,j}\}_{ i\ge
1,i+j\ge 2}$, and the multiplication in $G$ is defined by the formulas
\begin{equation*}
 e_{i,k}e_{j,l}=-{1\over2}e_{i+j,k+l-2} , \quad 
 f_{i,k}e_{j,l}=-{1\over2}f_{i+j,k+l-2}. \quad
 f_{i,k}f_{j,l}=0
\end{equation*}
If we let $t_i=e_{1,i}, u_i=f_{1,i}$ 
then $G$ is the algebra 
generated by $\{t_i,u_i\}_{i\ge 1}$, 
with the relations  $u_iu_j=0, t_i t_j = t_k t_l$ if $i+j=k+l$, and
$t_i u_j = t_k u_l$ if $i+j=k+l$.
\section{Algebraic Deformation Theory}
The discussion below holds for both \linf\ and \ainf\ algebras and the term
{\it infinity algebra} will be used to refer to either one.
Let $S$ be a graded commutative \k-algebra 
with identity,  containing  a distinguished maximal ideal 
${\frak m}\subset S$ such that $S/{\frak m}\cong\k$. 
Let $\varepsilon:S\ra S/{\frak m}=\k$ be the projection 
with $\varepsilon(1)=1$.   Let $V$ be an infinity algebra 
over the field \k, with codifferential $d$.
Then $V\tns S$ is an $S$-module and a deformation of $V$ with base
$S$ is defined to be an infinity algebra
structure over $S$ on  $V\tns S$, 
with the condition that 
$1\tns \varepsilon \colon V\tns S \to V\tns \k = V$ is an
infinity algebra homomorphism.

Suppose that $S$ is finite-dimensional and we have  a series  
$\tau = \sum\tau_n$  of  $S$-linear maps 
$\tau_n \colon \bigwedge^n(V\tns S )\to V\tns S$ which satisfy the
condition that $1\tns\varepsilon$ is an infinity algebra
homomorphism, that is 
\begin{equation}
(1\tns \varepsilon)\tau_n(v_1\tns s_1, \cdots, v_n\tns s_n)=
d_n(v_1,\cdots,v_n)\varepsilon (s_1\cdots s_n)
\end{equation}
for any
$v_1,\cdots,v_n \in V, s_1,\cdots,s_n \in S.$  Here, by homomorphism,
we are not implying that $\tau$ actually satisfies the codifferential
condition.
Then a linear map
$\alpha : \m^\ast \to C(V)$ can be defined by
\begin{equation}
\alpha(\ph)(v_1, \cdots, v_k)=(1\tns \ph)(\tau_k (v_1\tns 1, \cdots,
v_k\tns 1)- d_k(v_1,\cdots,v_k))
\end{equation}
for any linear functional $\ph\in \m^\ast$.  With this definition 
$\tau$ and $\alpha$
determine each other.  Let $F=\m^\ast,$ $\Delta$ be the
comultiplication in $F$ dual to the multiplication in \m, and 
$F_0=(\m/\m^2)^\ast.$
$F_0$ may be considered naturally as a subset of $F$ and in this sense 
$F_0\subset \ker \Delta.$  We will show that if $\tau$ is a codifferential,
providing $V\Tns S$ with structure of an infinity algebra,
 then $\delta\of\alpha
  +{1\over2}\mu\of(\alpha\tns\alpha)\of\Delta=0$,
so the  map
\begin{equation}
a\colon F_0\,
{\buildrel\alpha\over\longrightarrow}\, \ker \delta\,
{\buildrel{\pi}\over\longrightarrow}\, H(V)
\end{equation}
is well-defined.  In this
case, $a$ is called the {\it differential} of 
the deformation $\tau$.  An {\it infinitesimal
deformation of \g\ with base $S$} is defined to be any linear map 
$F_0 \to H(V)$. We will use proposition \ref{prp} below to show 
\begin{thm}
An infinitesimal deformation $a\colon F_0\to 
H(V)$ is the differential of some deformation with 
base $S$ if and only if $-{1\over2}a$ satisfies the condition of triviality of 
Massey $F$-products. 
\end{thm}
\begin{prp}\label{prp}
The operator $\tau$ satisfies the infinity algebra structure equation  
 if and only if $\alpha$ satisfies the Maurer-Cartan equation
 $$\delta\of\alpha
  +{1\over2}\mu\of(\alpha\tns\alpha)\of\Delta=0.$$
\end{prp}
\begin{prf} (\linf Algebra Case)
\label{MCprf2}
Let $\{m_i\}$ be a (homogeneous) basis of \m\  and let $\{m^i\}$ be the
dual basis of F.  
Suppose $m_i m_j ={\sum_{k}{c_{ij}^k m_k}}$. 
Then the comultiplication is given by $\Delta m^k ={\sum_{i,j}{(-1)^{m_i m_j}
 c_{ij}^k(m^i \tns m^j )}}$. 
Also for some 
 $\beta_n^i\in C^n(V)$, we have
$\tau_n(v_1\tns 1, \cdots, v_n\tns 1) = l_n(v_1,\cdots ,v_n)
 +\sum_i\beta_n^i(v_1,\cdots, v_n)\tns m_i$.
Thus for $v_1,\cdots,v_n \in V$, and $n=k+l-1$,
\allowdisplaybreaks
\begin{align*}
 \tau_k( \tau_l(&v_1,\cdots,v_l),v_{l+1},\cdots,v_n  )  \\
=&\tau_k \big( l_l(v_1,\cdots,v_l)
                   +\sum_i\beta_l^i(v_1,\cdots ,v_l)\tns m_i,
                   v_{l+1},\cdots,v_n \big) \\
=&\tau_k(l_l(v_1,\cdots,v_l),v_{l+1},\cdots,v_n) \\
 &   +\sum_i(-1)^{m_i(v_{l+1}+\cdots+v_n)}
       \tau_k(\beta_l^i(v_1,\cdots,v_l),v_{l+1},\cdots,v_n)\tns m_i \\
=&\l_k(l_l(v_1,\cdots,v_l),v_{l+1},\cdots,v_n) \\
 & +\sum_i\beta_k^i(l_l(v_1,\cdots,v_l),v_{l+1},\cdots,v_n)\tns m_i \\
 &+\sum_i(-1)^{m_i(v_{l+1}+\cdots +v_n)}
     l_k(\beta_l^i(v_1,\cdots,v_l),v_{l+1},\cdots,v_n)\tns m_i \\
 &+\sum_{r,s}(-1)^{m_s(v_{l+1}+\cdots+v_n)}
     \beta_k^r(\beta_l^s(v_1,\cdots,v_l),v_{l+1},\cdots,v_n)
     \tns m_rm_s \\
=&\l_k(l_l(v_1,\cdots,v_l),v_{l+1},\cdots,v_n) \\
 &+\sum_i\beta_k^i(l_l(v_1,\cdots,v_l),v_{l+1},\cdots,v_n)\tns m_i \\
 &+\sum_i(-1)^{m_i(v_{l+1}+\cdots +v_n)}
     l_k(\beta_l^i(v_1,\cdots,v_l),v_{l+1},\cdots,v_n)\tns m_i \\
 &+\sum_{i,r,s}(-1)^{m_s(v_{l+1}+\cdots v_n)}c_{rs}^i
     \beta_k^r(\beta_l^s(v_1,\cdots,v_l),v_{l+1},\cdots,v_n)
     \tns m_i
\end{align*}

If we use the notation $\alpha_k(m^i)=\alpha_k^i$ 
then 
\begin{equation*}
\beta_k^i(v_1,\cdots,v_k) = 
  (-1)^{m_i(k+1+v_1+\cdots+v_k)}
  \alpha_k^i(v_1,\cdots,v_k).
\end{equation*}
Substituting this into the previous
expression yields
\begin{align*}
  \l_k(l_l(v_1,\cdots,v_l),v_{l+1},\cdots,v_n) 
  +\sum_i (-1)^{m_i(n+v_1+\cdots+v_n)} 
  M_i(v_1,\cdots,v_n) \tns m_i,
\end{align*}
where 
\begin{multline*}
  M_i(v_1,\cdots,v_n) = 
     \alpha^i_k(l_l(v_1,\cdots,v_l),v_{l+1},\cdots,v_n)\\
   +(-1)^{km_i}
     l_k(\alpha^i_l(v_1,\cdots,v_l),v_{l+1},\cdots,v_n)\\
   +\sum_{r,s}c_{rs}^i(-1)^{m_rm_s + km_s}
     \alpha_k^r(\alpha_l^s(v_1,\cdots,v_l),v_{l+1},\cdots,v_n). 
\end{multline*}
The requirement for $\tau$ to be a codifferential is 
that $\{ \tau,\tau \}=0$ or equivalently, 
 ${\displaystyle \sum_{k+l=n+1} \{ \tau_k,\tau_l \} = 0}$,
for all $n$.  
Since
$(-1)^{kl+(k+1)(l+1) +(k-1)l} = (-1)^{(l-1)k}$
we have
\begin{multline*}
\{\tau_k,\tau_l \}(v_1,\cdots,v_n) = \\
\sum_{\sigma\in\sh(l,k-1)} \sigsign (-1)^{(k-1)l}
  \tau_k(\tau_l(\permarg{v}{1}{l}),\permarg{v}{l+1}{n})\\
+\sum_{\sigma\in\sh(k,l-1)} \sigsign (-1)^{(l-1)k}
  \tau_l(\tau_k(\permarg{v}{1}{k}),\permarg{v}{k+1}{n}),
\end{multline*}
so
\begin{multline}
\label{reduction}
\{\tau_k, \tau_l\} + \{\tau_l,\tau_k\} \\
 = 2 \sum_{\sigma\in\sh(l,k-1)} \sigsign (-1)^{(k-1)l}
  \tau_k(\tau_l(\permarg{v}{1}{l}),\permarg{v}{l+1}{n}).
\end{multline}
Because of this fact, $\tau$ gives an \linf\ algebra structure exactly when
\begin{multline*}
\sum_{k+l=n+1 \atop \sigma \in \sh(l,k-1)}
 \sigsign (-1)^{(k-1)l} 
  \Big(  \l_k(l_l(v_1,\cdots,v_l),v_{l+1},\cdots,v_n) \\
 \left.  +\sum_i (-1)^{m_i(n+v_1+\cdots+v_n)} 
  M_i(v_1,\cdots,v_n) \tns m_i \right)=0.
\end{multline*}
Since the first term is zero by the generalized Jacobi identity and 
 $\{m^i\}$ is a basis, this is the same as requiring 
\begin{equation*}
\sum_{k+l=n+1 \atop \sigma \in \sh(l,k-1)}
(-1)^{(k-1)l}\sigsign M_i(\permarg{v}{1}{n}) = 0, 
\end{equation*}
for all $i$ and $n$.

Looking at the terms in the Maurer-Cartan equation, 
\begin{equation*}
\delta\of\alpha(m^i)=\delta\of\alpha^i= \{ l, \alpha^i  \}
 =\sum_{k,l} \{ l_k,\alpha^i_l \}.
\end{equation*}
Thus
\begin{multline*}
\delta\of\alpha(m^i)(v_1,\cdots,v_n)
= \sum_{k+l=n+1} \{ l_k,\alpha^i_l  \} (v_1,\cdots,v_n)\\
= \sum_{k+l=n+1 \atop \sigma\in\sh(l,k-1)}
   \sigsign(-1)^{(k-1)(l+m_i)}
  l_k\big(\alpha^i_l(\permarg{v}{1}{l}),\permarg{v}{l+1}{n}\big)\\
 + \sum_{k+l=n+1 \atop\sigma\in\sh(k,l-1)}
  \sigsign(-1)^{(l+1)k+m_i}
  \alpha^i_l\big(l_k(\permarg{v}{1}{k}),\permarg{v}{k+1}{n}\big)
\end{multline*}
and 
\begin{align*}
{1\over2}\mu\of(\alpha\tns\alpha)\of\Delta(m^i)
&={1\over2}\sum_{r,s}(-1)^{m_rm_s}c^i_{rs}
   \mu\of(\alpha\tns\alpha)(m^r\tns m^s)\\
&={1\over2}\sum_{r,s,k,l}(-1)^{m_rm_s+m_r}c^i_{rs}
  \mu\of(\alpha^r_k\tns\alpha^s_l)\\
&={1\over2}\sum_{r,s,k,l}(-1)^{m_r(m_s+1)}c^i_{rs}
  \{ \alpha^r_k,\alpha^s_l \}
\end{align*}
so, by an identity similar to (\ref{reduction}),
\begin{align*}
{1\over2}\mu\of(\alpha\tns\alpha)\of\Delta(m^i)
  (v_1,\cdots,v_n)\\
={1\over2}\sum_{r,s \atop k+l=n+1}(-1)^{m_r(m_s+1)}&c^i_{rs}
  \{ \alpha^r_k,\alpha^s_l \}(v_1,\cdots,v_n)\\
=\sum_{{r,s\atop k+l=n+1}\atop\sigma\in\sh(l,k-1)}
 \sigsign & c^i_{rs} (-1)^{m_r(m_s+1)+(k-1)(l+m_s)} \\
  \times \alpha^r_k(\alpha^s_l&
    (\permarg{v}{1}{l}),\permarg{v}{l+1}{n})
\end{align*}
Putting the pieces together,
\begin{multline*}
(\delta\of\alpha+{1\over2}\mu\of(\alpha\tns\alpha)
 \of\Delta)(m^i)(v_1,\cdots,v_n)\\
=\sum_{\sigma\in\sh(l,k-1)}\sigsign 
 M_i(\permarg{v}{1}{n}).
\end{multline*}
Hence the proposition is proved for \linf\ algebras.
\end{prf}
\begin{prf} (\ainf\ Algebra Case)
\label{MCprf3}
Let $\{m_i\}$ be a (homogeneous) basis of \m\  and let $\{m^i\}$ be the
dual basis of F.  
Suppose that $m_i m_j = \displaystyle{\sum_{k}{c_{ij}^k m_k}}$. 
Then $\Delta m^k =  \displaystyle{\sum_{i,j}{(-1)^{m_i m_j}
 c_{ij}^k(m^i \tns m^j )}}$. 
Also for $v_1,\cdots,v_n \in V$, 
$$\tau_n(v_1\tns 1, \cdots, v_n\tns 1) = \mu_n(v_1,\cdots ,v_n)
 +\sum_i\beta_n^i(v_1,\cdots, v_n)\tns m_i,$$
for some 
 $\beta_n^i\in C^n(V).$

Thus for $v_1,\cdots,v_n \in V$, any $j$, and $k+l-1=n$,
\begin{align*}
\tau_k(&v_1,\cdots,v_{j-1},\tau_l(\midv),\rightv)\\
=&\tau_k\big( \leftv, \mu_l(\midv)
                  +\sum_i\beta_l^i(\midv)\tns m_i,\rightv \big)\\ 
=&\mu_k(\leftv,\mu_l(\midv),\rightv)\\  
&+\sum_i\beta^i_k(\leftv,\mu_l(\midv),\rightv)\otimes m_i\\
&+\sum_i  
\s{m_ix}
 \mu_k(\leftv,\beta^i_l(\midv),\rightv)\otimes m_i\\
 &+\sum_{r,s}
\s{m_sx}
    \beta^r_k(\leftv,\beta^s_l(\midv),\rightv)\otimes m_rm_s
\\
=&\mu_k(\leftv,\mu_l(\midv),\rightv)\\  
&+\sum_i\beta^i_k(\leftv,\mu_l(\midv),\rightv)\otimes m_i\\
&+\sum_i 
\s{m_i x}
 \mu_k(\leftv,\beta^i_l(\midv),\rightv)\otimes m_i\\
 &+\sum_{r,s}
\s{m_sx}
c_{rs}^i
    \beta^r_k(\leftv,\beta^s_l(\midv),\rightv)\otimes m_i
\end{align*}
where $x=v_{j+l+1}+\cdots+v_n$.

As before, if we use the notation $\alpha_k(m^i)=\alpha_k^i$ 
then 
\begin{equation*}
\beta_k^i(v_1,\cdots,v_k) = 
  (-1)^{m_i(k+1+v_1+\cdots+v_k)}
  \alpha_k^i(v_1,\cdots,v_k)
\end{equation*}
and substituting this into the previous
expression yields
\begin{multline*}
  \mu_k(\leftv,\mu_l(\midv),\rightv)\\ 
  +\sum_i (-1)^{m_i(n+v_1+\cdots+v_n)} 
  M_i(v_1,\cdots,v_n) \tns m_i,
\end{multline*}
where 
\begin{align*}
  M_i(&v_1,\cdots,v_n) \\ 
 =&\alpha^i_k(\leftv,\mu_l(\midv),\rightv)\\
  &+(-1)^{m_i(k+v_1+\cdots +v_{j-1})}
  \mu_k(\leftv,\alpha^i_l(\midv),\rightv)\\
  &+\sum_{r,s}c_{rs}^i 
\s{x}
   \alpha^r_k(\leftv,\alpha^s_l(\midv),\rightv),
\end{align*}
where $x=m_s(m_r +k+v_1+\cdots+v_{j-1})$.
The requirement for $\tau$ to be a codifferential is
that $\{ \tau,\tau \}=0$ or equivalently, 
${\displaystyle \sum_{k+l=n+1} \{ \tau_k,\tau_l \} = 0}$,
for all $n$.  By reasoning similar to (\ref{reduction}) this may be
expressed as
\begin{equation*}
\sum_{1\leq j\leq k }
\s{x}
   \tau_k(\leftv,\tau_l(\midv),\rightv)=0,
\end{equation*}
where $x=l(k+v_1+\cdots+v_{j-1}) + j(l-1)$,
or,
\begin{multline*}
\sum_{ k+l=n+1\atop 1\leq j\leq k }
\s{x}
 \Big( (\mu_k(\leftv,\mu_l(\midv),\rightv)  \\
 \left. + \sum_i (-1)^{m_i(n+v_1+\cdots+v_n)} 
 M_i(v_1,\cdots ,v_n) \tns m_i \right)=0,
\end{multline*}
where $x=l(k+v_1+\cdots+v_{j-1}) + j(l-1)$.
The first term is zero because of the \ainf\ algebra structure on $V$ and
since $\{m^i\}$ is a basis, this may be written as,
\begin{equation*}
\sum_{k+l=n+1 \atop 1\leq j\leq k}
 (-1)^{(k-1)l + l(v_1+\cdots +v_{j-1})+(j-1)(l-1)}
 M_i(v_1,\cdots ,v_n)=0,
\end{equation*}
for all $i$  and $n$.

From the Maurer-Cartan equation we have
\begin{align*}
\delta\of\alpha(m^i)(v_1, &\cdots ,v_n)\\
&= \sum_{k+l=n+1} \{ \mu_k,\alpha^i_l  \} (v_1,\cdots ,v_n)\\
&=\sum_{k+l=n+1 \atop 1\leq j \leq k}
\s{x}
\left( \mu_k(\leftv,\alpha^i_l(\midv),\rightv)\right. \\
&{\hskip .25in}+
\s{y}
 \left. \alpha^i_k(\leftv,\mu_l(\midv),\rightv)\right),
\end{align*}
 where $x=(l+m_i)(k+v_1+\cdots+v_{j-1}) + j(l-1)+m_i+1$,
$y =m_i(k + v_1 +\cdots +v_{j-1})$,
and
\begin{align*}
{1\over2}\nu\of(\alpha\tns\alpha)\of\Delta(m^i)
  (v_1,\cdots,&v_n)\\
={1\over2}\sum_{r,s \atop k+l=n+1}(&-1)^{m_r(m_s+1)}c^i_{rs}
  \{ \alpha^r_k,\alpha^s_l \}(v_1,\cdots,v_n)\\
=\sum_{{r,s \atop k+l=n+1} \atop 1\leq j \leq k}c^i&_{rs}
(-1)^{m_rm_s+k(l+m_r)+(l+m_s)(v_1+\cdots v_{j-1}) +j(l-1)} \\
&\times \alpha^r_k(\leftv,\alpha^s_l(\midv),\rightv).
\end{align*}
Together we have,
\begin{multline*}
(\delta\of\alpha+{1\over2}\nu\of(\alpha\tns\alpha)
 \of\Delta)(m^i)(v_1,\cdots,v_n)\\
=\sum_{ k+l=n+1 \atop 1\leq j\leq k }
 (-1)^{(k-1)l + l(v_1+\cdots +v_{j-1})+(j-1)(l-1)}
 M_i(v_1,\cdots ,v_n).
\end{multline*}
Thus the proposition is proved for \ainf\ algebras.
\end{prf}
\section{Conclusion}
When the examples of Massey F-products 
in \cite{fl} were gathered together, the examples which we present here
were known to the authors, but it was decided that the extra definitions
necessary to appreciate the deformation theory of infinity algebras
merited a separate treatment.  At the same time, we did not know how
to extend the notion of deformation of algebras with a base to infinity
algebras. When we discovered how to do this, we decided that the combination
of results might be interesting, and so we have collected them in this
article. We would like to take this opportunity to thank Dmitry Fuchs,
who encouraged us to write this paper, and who offered some insightful
remarks.
\providecommand{\bysame}{\leavevmode\hbox to3em{\hrulefill}\thinspace}


\begin{thebibliography}{10}

\bibitem{aksz}
M.~Alexandrov, M.~Kontsevich, A.~Schwarz, and O.~Zaboronsky, \emph{The geometry
  of the master equation and topological quantum field theory}, Preprint, 1995.

\bibitem{ber}
F.A. Berends, G.J.H. Burgers, and H.~van Dam, \emph{On the theoretical problems
  in constructing interactions involving higher-spin massless particles},
  Nuclear Physics B \textbf{260} (1985), no.~2, 295--322.

\bibitem{conn}
A.~Connes, \emph{Non-commutative differential geometry}, Institute Des Hautes
  Etudes Scientifiques \textbf{62} (1985), 41--93.

\bibitem{fl}
D.~Fuchs and L.~Lang, \emph{{M}assey products and deformations}, Preprint
  q-alg/9602024, 1996.

\bibitem{getz}
E.~Getzler and J.D.S. Jones, \emph{\hbox{$A_\infty$}-algebras and the cyclic
  bar complex}, Illinois Journal of Mathematics \textbf{34} (1990), no.~2,
  256--283.

\bibitem{getz2}
\bysame, \emph{Operads, homotopy algebra, and iterated integrals for double
  loop spaces}, Preprint, 1995.

\bibitem{hoch}
G.~Hochschild, \emph{On the cohomology groups of an associative algebra},
  Annals of Mathematics \textbf{46} (1945), 58--67.

\bibitem{kast}
D.~Kastler, \emph{Cyclic cohomology within the differential envelope}, Travaux
  En Cours, Hermann, Paris, 1988.

\bibitem{lm}
T.~Lada and M.~Markl, \emph{Strongly homotopy lie algebras}, Comm. in Algebra
  \textbf{23} (1995), 2147--2161.

\bibitem{ls}
T.~Lada and J.~Stasheff, \emph{Introduction to sh lie algebras for physicists},
  Preprint hep-th/9209099, 1990.

\bibitem{lod}
J.~Loday, \emph{Cyclic homology}, Springer-Verlag, 1992.

\bibitem{mar}
M.~Markl, \emph{A cohomology theory for {A}($m$)-algebras and applications},
  Journal of Pure and Applied Algebra \textbf{83} (1992), no.~6, 141--175.

\bibitem{pen1}
M.~Penkava, \emph{\hbox{$L_\infty$} algebras and their cohomology}, Preprint
  q-alg/9512014, 1995.

\bibitem{ss}
M.~Schlessinger and J.~Stasheff, \emph{The {L}ie algebra structure of tangent
  cohomology and deformation theory}, Journal of Pure and Applied Algebra
  \textbf{38} (1985), 313--322.

\bibitem{seib}
P.~Seibt, \emph{Cyclic homology of algebras}, World Scientific, 1987.

\bibitem{sta1}
J.D. Stasheff, \emph{On the homotopy associativity of {H}-spaces {I}},
  Transactions of the AMS \textbf{108} (1963), 275--292.

\bibitem{sta2}
\bysame, \emph{On the homotopy associativity of {H}-spaces {II}}, Transactions
  of the AMS \textbf{108} (1963), 293--312.

\bibitem{sta3}
J.D. Stasheff, \emph{Closed string field theory, strong homotopy lie algebras
  and the operad actions of moduli spaces}, Preprint hep-th/9304061, 1992.

\bibitem{umb}
R.~Umble, \emph{The deformation complex for differential graded hopf algebras},
  Preprint, 1994.

\bibitem{WitZwie}
E.~Witten and B.~Zwiebach, \emph{Algebraic structures and differential geometry
  in two-dimensional string theory}, Nuclear Physics B \textbf{377} (1992),
  55--112.

\bibitem{Z}
B.~Zwiebach, \emph{Closed string field theory: {Q}uantum action and the
  {B}atalin-{V}ilkovisky master equation}, Nuclear Physics B \textbf{390}
  (1993), 33--152.

\end{thebibliography}
\end{document}